\documentclass[twocolumn]{autart}  

\usepackage{circuitikz}

\usepackage{setspace}
\usepackage{graphicx}
\usepackage{graphics}
\usepackage{amsfonts}
\usepackage{epsfig}
\usepackage{mathrsfs}
\usepackage{amsmath}
\usepackage{amssymb}
\usepackage{color}
\usepackage{bbold}
\usepackage{dsfont}
\usepackage{siunitx}
\usepackage{flushend} 

\def\qedp{\hspace*{\fill}~{\tiny $\blacksquare$}}

\DeclareSymbolFont{bbold}{U}{bbold}{m}{n}
\DeclareSymbolFontAlphabet{\mathbbold}{bbold}

\newtheorem{example}{Example} 

\newtheorem{remark}{Remark} 

\newcommand{\be}{\ensuremath{\beta}}

\usepackage[prependcaption,colorinlistoftodos]{todonotes}

\newcommand{\blue}[1]{\textcolor{blue}{#1}}

\def\be{\begin{equation}}
\def\ee{\end{equation}}
\def\ba{\begin{array}}
\def\ea{\end{array}}
\def\eqa{\begin{eqnarray}}
\def\eqe{\end{eqnarray}}

\begin{document}

\begin{frontmatter}
\title{A power consensus algorithm for DC microgrids\thanksref{footnoteinfo}}

\thanks[footnoteinfo]{An abridged version of the paper has been submitted to 
the 20th IFAC World Congress, 9-14 July 2017, Toulouse, France.
}
\author[groningen]{C.~De Persis}\ead{c.de.persis@rug.nl},  
\author[groningen]{E.R.A.~Weitenberg}\ead{e.r.a.weitenberg@rug.nl},
\author[zurich]{F.~D\"{o}rfler}\ead{dorfler@ethz.ch}

\address[groningen]{
Engineering and Technology Institute Groningen, 
University of Groningen,
Nijenborgh 4, 
9747 AG Groningen,
The Netherlands
}

\address[zurich]{
Automatic Control
Laboratory, Swiss Federal Institute of Technology, 8092
Zurich, Switzerland}


\begin{abstract}    
A novel power consensus algorithm for DC microgrids is proposed and analyzed. DC microgrids are networks {composed of DC sources, loads, and interconnecting lines.}  They are represented by differential-algebraic equations connected over an undirected weighted graph that models the electrical circuit. A second graph represents the communication network  over which the source nodes  exchange information about the instantaneous powers, which is used to adjust the injected current accordingly. This give rise to a nonlinear consensus-like system of differential-algebraic equations that is analyzed via Lyapunov functions inspired by the physics of the system. We establish convergence to the set of equilibria consisting of weighted consensus power vectors as well as preservation of the {weighted} geometric mean of the source voltages. The results apply to networks with constant impedance, constant current and constant power loads. 
\end{abstract}


\end{frontmatter}

\section{Introduction}
The proliferation  of renewable energy sources and storage devices  that are intrinsically operating using the DC regime is stimulating interest in the design and operation of DC microgrids, which have the additional desirable feature of  preventing the use of inefficient power conversions at different stages. These DC microgrids might have to be deployed in areas where an AC microgrid is already in place, creating what is called a hybrid microgrid \cite{loh}, for which rigorous analytical studies are still in their infancy. 
 Furthermore, the envisioned future in which power generation is far away from the major consumption sites raises the problem of how {to transmit} power with low losses, a problem for which High Voltage Direct Current (HVDC) networks perform comparatively better than AC networks.  Finally, also mobile grids on
ships, aircrafts, and trains are based on a DC architecture.

With DC and hybrid microgrids, as well as  HVDC networks, on the rise, 
we need to develop a deeper system-theoretic understanding {of this interesting class of dynamical networks}. 
In this paper we propose and analyse a control algorithm for a DC microgrid that enforces power sharing among the different power sources. 

\subsection{Literature review}
The literature on DC microgrids is rapidly growing. We summarize below the  contributions that share a systems and control{-}theoretic point of view on these networks. 
The work \cite{hill2} relies on a  cooperative control paradigm for {DC} microgrids to replace the conventional secondary control by a voltage and a current regulator.
In \cite{zhao.dorfler.aut15} a voltage droop controller for DC microgrids inspired by frequency droop in AC power networks is analyzed\blue{,} and a secondary consensus control strategy is added to prevent voltage drift and achieve optimal current injection. The paper \cite{belk2016stability} models the DC microgrid via the Brayton-Moser equations and uses this formalism to show that with the addition of a decentralized integral controller voltage regulation to a desired reference value is achieved. Other schemes achieving desirable power sharing properties are proposed but no formal analysis is provided. In \cite{tucci2016consensus}, a  secondary consensus-based control scheme for current sharing and voltage balancing in DC microgrids is designed  in a Plug-and-Play fashion to allow for the addition or removal of generation units.  A distributed control method to enforce  power sharing among a cluster of {DC} microgrids is proposed in \cite{moayedi.davoudi.2016}. 
Other work has  focused on the challenges in the stability analysis of DC microgrids using consensus-like algorithms due to the interaction between the communication network and the physical one \cite{meng.tsg2016}. Finally, feasibility of the nonlinear algebraic equations in DC power circuits is studied by \cite{barabanov}, \cite{john.resistive.nets}, and \cite{lavei.rantzer}.

A closely related research area is that of multi-terminal HVDC transmission systems. The paper \cite{sarlette.aut12} focuses on cooperative frequency control for these networks.  In \cite{andreasson.hvdc} distributed controllers that keep the voltages close to a nominal value and guarantee a fair power sharing are considered, whereas passivity-based decentralized PI control for  the global asymptotic stabilisation of multi-terminal high-voltage is studied in \cite{zonetti.cep2015}. 
The paper \cite{DZ-RO-JS:16} studies feasibility and power sharing under decentralized droop control.
We refer to \cite[Chapter 4]{Zonetti.thesis} for an annotated bibliography of HVDC transmission systems. 

\subsection{Main contribution} {This} paper focuses on a new control algorithm to stabilize a DC microgrid under different
load characteristics
while achieving power sharing among the sources. Our controller is enabled by  communicating the instantaneous source power measurements among neighboring source nodes, averaging these measurements and setting the voltage at the source terminals accordingly. An additional feature of the algorithm is that {a weigthed} geometric average of the source  voltages is preserved. 

The system dynamics present interesting features. {By} averaging the power measurements that  the sources communicate amongst each other, the system dynamics becomes an intriguing combination of the physical network (the weighted Laplacian of the electrical circuit appearing in the power measurements) and the communication network (over which the information about the power measurements is exchanged). ``ZIP" (constant impedance, constant current and constant power) loads introduce algebraic equations in the system's dynamics, adding additional complexity and nonlinearities.

To analyze this  system of nonlinear differential-algebraic equations without going through a linearization of the dynamics, {Lyapunov-based} arguments become very convenient. The Lyapunov functions in this case are constructed starting from the power dissipated in the network that is further shaped to take into account the specifics of the dynamics. The presence of the loads, which shift the equilibrium of interest, is taken into account by the so-called Bregman  function \cite{claudio.nima.ecc16}.  The level sets of the Lyapunov functions are used to estimate the excursion of the state response of these systems and therefore, combined with the preservation of the geometric average of the source voltages, can be used to obtain an estimate of the voltage at steady state. 

{Reactive} power sharing algorithms have been first suggested by \cite{schiffer.tcst16} for network-reduced AC microgrids whose voltage dynamics show similar features as in DC grids. In this paper we show that a similar idea can be adopted also for  network preserved DC microgrids. The novelties of this contribution with respect to \cite{schiffer.tcst16} are the different dynamics of the system under study, the explicit consideration of algebraic equations in the model and the use of Lyapunov arguments to prove the main results.

\subsection{Paper organization} The model of the DC microgrid is introduced in Section \ref{sec.dc.microgrid}. The power consensus algorithm is introduced in Section \ref{sec.power.controller}. The analysis of the closed-loop system is  carried out in Section \ref{sec.power.consensus.constant.current.load} 
{for the general case of ZIP loads, and then specialised to the case of ZI loads, since the latter  permits to obtain stronger results under weaker conditions. }
{The simulations of the algorithm are provided} in Section \ref{simulations}. 
Conclusions are drawn in Section \ref{sec.conclusions}. 

\subsection{Notation} Given a vector $v$, the symbol $[v]$ represents the diagonal matrix whose diagonal entries are the components of $v$. The notation ${\rm col}(v_1, v_2, \ldots, v_n)$, with $v_i$ scalars,  represents the vector $[\begin{matrix} v_1 & v_2 & \ldots & v_n\end{matrix}]^T$. If $v_i$ are matrices having the same number of columns, then ${\rm col}(v_1, v_2, \ldots, v_n)$ denotes the matrix 
$[\begin{matrix} v_1^T & v_2^T & \ldots & v_n^T\end{matrix}]^T$. The symbol $\mathbb{1}_n$ represents the $n$-dimensional  vector of all $1$'s, whereas $\mathbb{0}_{m\times n}$ is the $m\times n$ matrix of all zeros. When the size of the matrix is clear from the context the index is omitted. The $n\times n$ identity matrix is represented as $\mathbb{I}_n$. Given a vector $v\in\mathbb{R}^n$, the symbol $\boldsymbol{\ln}(v)$ denotes the element-wise logarithm, i.e., the vector $[\,\ln(v_1)\ldots\ln(v_n)\,]^T$.  

\section{DC resistive microgrid}\label{sec.dc.microgrid}
The DC microgrid is modeled as an undirected connected graph $\mathcal{G}=(\mathcal{V},\mathcal{E})$, with $\mathcal{V}:=\{1,2,\ldots,n\}$ the set of nodes {(or buses)} and 
$\mathcal{E}\subseteq \mathcal{V}\times\mathcal{V}$ the set of edges. The edges represent the {interconnecting} lines of the microgrid, which we assume here to be resistive. Associated to each edge is a weight modeling {the}
conductance (or reciprocal resistance) $1/r_k>0$,
with $k\in \mathcal{E}$.
The set of nodes is partitioned into the two subsets of $n_s$ {DC} sources $\mathcal{V}_s$ and $n_l$ loads $\mathcal{V}_l$, with $n_s+n_l=n$.  

The current-potential relation in a resistive network is given by the identity $I= B\Gamma B^T V$, with {$B \in \mathbb R^{n \times |\mathcal E|}$ being} the incidence matrix of $\mathcal{G}$ and $\Gamma={\rm diag}\{r_1^{-1}, \ldots,r_m^{-1}\}$ the diagonal {matrix of conductances}. Considering the partition of the nodes in sources and loads, the relation  rewrites as 
\be\label{resistive.network}
\begin{bmatrix}
I_s \\ I_l
\end{bmatrix}
=
\begin{bmatrix}
B_s \Gamma B_s^T & B_s \Gamma B_l^T\\
B_l \Gamma B_s^T & B_l \Gamma B_l^T\\
\end{bmatrix}
\begin{bmatrix}
V_s \\ V_l
\end{bmatrix}
=: 
\begin{bmatrix}
Y_{ss} & Y_{sl}\\
Y_{ls} & Y_{ll}
\end{bmatrix}
\begin{bmatrix}
V_s \\ V_l
\end{bmatrix},
\ee
where $I_s={\rm col} (I_{1},\dots,I_{n_s})$, $I_l={\rm col} (I_{n_s+1},\dots,I_{n})$, $V_s={\rm col} (V_{1},\dots,V_{n_s})$, $V_l={\rm col} (V_{n_s+1},\dots,V_{n})$ and $B= {\rm col}(B_s, B_l)$.

{Observe that both $Y_{ss}$ and $Y_{ll}$ are positive definite since they are principal submatrices of a Laplacian of a connected undirected graph. This allows us to eliminate the load voltages as $V_{l} = Y_{ll}^{-1}I_{l} - {Y_{ll}^{-1} Y_{ls} } V_{s}$ and reduce the network to the source nodes $\mathcal V_{s}$ with balance equations
\begin{equation}
I_{s} - {Y_{ll}^{-1} Y_{ls} } I_{l} ={Y_{red}} V_{s}
\label{eq:Kron-reduction}
\,,
\end{equation}
where  ${Y_{red}}  = Y_{ss}- {Y_{ll}^{-1} Y_{ls} } Y_{ls}$ is known as the Kron-reduced conductance matrix \cite{FD-FB:11d} and $- {Y_{ll}^{-1} Y_{ls} } I_{l}$ is the mapping of the load current injections to the sources.
}

\section{Power consensus controllers}\label{sec.power.controller}
We propose controllers that force the different sources to share the total power injection {in prescribed ratios} \cite{schiffer.tcst16}. For this purpose,
a communication network is deployed to connect the source nodes,  through  which the controllers exchange information about the instantaneous injected powers. This communication network is modelled as an undirected unweighted graph $(\mathcal{V}_c, \mathcal{E}_c)$, where $\mathcal{V}_c=\mathcal{V}_s$. Associated with the communication graph is the $n_s\times n_s$ Laplacian matrix $L_c=D_c-A_c$, where $D_c$ is the degree matrix and $A_c$ is the adjacency matrix of the communication graph. Note that the nodes of the communication network (but not necessarily the edges) coincide with the source nodes of the microgrid. For each node $i\in \mathcal{V}_s$, the set $\mathcal{N}_{c,i}=\{j\in \mathcal{V}_s: \{i, j\}\in \mathcal{E}_c\}$ represents the neighbors connected to node $i$ via the communication graph.

{\it Controllers.} The proposed controllers are of the form
\be\label{circuit.interpretation.control}
\mathcal{C}_{i}(V_{i}) \dot V_{i}  = -I_{i}+ u_{i},\quad i\in \mathcal{V}_s,
\ee
where 
\be\label{nonlinear.capacitance}
\mathcal{C}_{i}(V_{i})= V_{i}^{-2} D_{ci}^{-1} C_{i}^2,\quad i\in \mathcal{V}_s
\ee 
can be interpreted as a nonlinear capacitance, 
$C_{i}>0$ 
is a  positive parameter {of suitable units such that $\mathcal{C}_{i}(V_i)$ actually has the units of a capacitance,} 
$I_{i}$ is the  injected current at node $i\in \mathcal{V}_s$ as defined in \eqref{resistive.network},  and the term
\be\label{u.s}
u_{i}= V_{i}^{-1} D_{ci}^{-1} C_{i} \sum_{j\in \mathcal{N}_{c,i}} C_{j}^{-1}
P_{j},\quad i\in \mathcal{V}_s
\ee 
represents an ideal current source that is controlled as a function of the local voltage $V_{i}$ and the injected power $P_{j}= V_{j} I_{j}$ at the neighboring node sources  $j\in \mathcal{N}_{c,i}$.

{The  dynamic controllers \eqref{circuit.interpretation.control}--\eqref{u.s} are initialised at positive values of the voltage, that is $V_i(0)>0$ for all $i\in \mathcal{V}_s$. It will be made evident in  later sections that these controllers {render} the positive orthant $\mathbb{R}^{n_s}_{>0}$ positively invariant, thus showing that the positivity of the initial source voltages yields positivity of these variables for all $t\ge 0$.}

\begin{remark} {\bf (Circuit interpretation)} 
The control algorithm has the circuit interpretation given in Fig.~\ref{fig.control}.
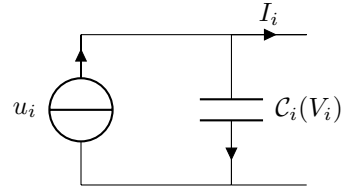
\begin{figure}[h!]
  \begin{center}
    \begin{circuitikz}
      \draw (0,0)
      to[I,i=$u_{i}$] (0,2) 
      to[short] (2,2)
      to[C=$\mathcal{C}_{i}(V_{i})$, i=$$] (2,0) 
      to[short] (0,0);
      \draw (2,2)
      to[short, i=$I_{i}$] (3,2);
      \draw (2,0)
        to (3,0);
    \end{circuitikz}
    \caption{A circuit interpretation of the controller \eqref{circuit.interpretation.control}. \label{fig.control}}
  \end{center}
\end{figure}
Comparing with \cite[(4)]{belk2016stability}, the ideal current source $u_{i}$ can be generated also by a voltage source with value  $v_{i}$ in series with a resistance $r_{i}$ provided that $v_{i} = r_{i} u_{i} + V_{i}$. Finally, the dynamic droop controller in \cite{zhao.dorfler.aut15} corresponds in our notation to a constant capacitance $\mathcal C_{i}$ and current source $u_{i}$.
\end{remark}
\begin{remark} {\bf (Digital implementation)}The control algorithm above need not be implemented analogically. In its digital implementation, information regarding the injected powers of neighbouring sources,
$P_{j}$, $j\in \mathcal{N}_{c,i}$, 
are broadcasted and processed along with the current measurement $I_{i}$ to compute the source voltage value $V_{i}$ applied at the source terminals. 
\end{remark}
Multiplying both sides of \eqref{circuit.interpretation.control} by $V_{i}^2D_{c,i}C_{i}^{-1}$, one arrives at the closed-loop system
\begin{eqnarray}
C_{i}\dot V_{i} &=& - V_{i} D_{ci}  C_{i}^{-1} P_{i} + V_{i}\sum_{j\in \mathcal{N}_{c,i}} C_{j}^{-1} P_{j}\nonumber\\
&=& 
V_{i}\sum_{j\in \mathcal{N}_{c,i}} ( C_{j}^{-1}  P_{j}- C_{i}^{-1}   P_{i}), \quad i\in \mathcal{V}_s,
\label{consensus.power.control}
\end{eqnarray}
that is, the voltage at the source terminal is updated according to a weighted power consensus algorithm scaled by the voltage. 
%
Provided that $V_i\ne 0$ (a property that will be established in the next sections), equation \eqref{consensus.power.control} shows that at steady state the algorithm achieves {\em proportional power sharing} according to the $C_i$ ratios, namely
\be\label{pps}
\frac{P_{j}}{C_{j}}= \frac{P_{i}}{C_i}, \quad \forall i,j\in \mathcal{V}_s.
\ee
A detailed characterisation of the steady-state power signals  is given in the next section (Lemma \ref{lem.equilibria}). 

For interpretation purposes, we write  \eqref{consensus.power.control}\,as
\begin{equation*}
	\frac{d}{dt} \, C_{i}\boldsymbol{\ln}(V_{i})
	= \sum_{j\in \mathcal{N}_{c,i}} ( C_{j}^{-1}  P_{j}- C_{i}^{-1}   P_{i}), \quad i\in \mathcal{V}_s\,.
\end{equation*}
In a classic power system analysis \cite{H-DC:11}, the term $C_{i}\boldsymbol{\ln}(V_{i})$ is the natural energy representation of a power source of {constant value} $C_{i}$. {The interpretation of the closed loop \eqref{consensus.power.control} is then that {the voltage at} this constant power source is adapted according to a power consensus algorithm.}

\begin{remark}{\bf ({Alternative} power sharing control)}\label{rem.dapi}
A possibly more simplistic and obvious power sharing controller {inspired by the current-sharing controller in \cite{zhao.dorfler.aut15} is} based on a distributed averaging integral control {given by}%
\be\ba{rl}\label{dapi}\
C_i \dot V_i &= -I_i + p_i \\
D_i \dot p_i &= I_i - p_i + \sum_{j\in\mathcal{N}_{c,i}} {(C_j^{-1} V_j p_j - C_i^{-1} V_i p_i)}, i\in \mathcal{V}_s
\ea\ee
where $p_{i}$ is a {control variable} in units of currents.
Any steady  state of this controller would guarantee for all $i \in \mathcal{V}_s$ that {$\dot V_{i} = 0$, and $p_{i}=I_{i}$ is the steady-state current injection,} and the vector of  power injections $C_s^{-1} [V_{s}] p$ has all identical entries (power sharing). Numerical results (see Section \ref{simulations}) show that \eqref{consensus.power.control} and \eqref{dapi} perform similarly. {Indeed, in the limit $D_{i}=0$, near steady-state, and for nearly unit voltages (in per unit system), the closed-loops  \eqref{dapi} and \eqref{consensus.power.control} have similar dynamics. In} the rest of the paper we focus on the analysis of \eqref{consensus.power.control}. 
\end{remark}

%
%
%
{\it Loads.} Depending on the particular load models, the term $I_l$ in \eqref{resistive.network} takes different expression and will henceforth be denoted as $I_l(V_l)$ to stress the functional dependence on the load voltages. 
{Prototypical} {load models that are of interest include the following:}
\begin{itemize}
\item[(i)] constant current loads: $I_l(V_l)= I_l^*\in \mathds{R}_{<0}^{n_l}$,
\item[(ii)] constant impedance: $I_l(V_l)= {-Y_l^{*}} V_l$, with $Y_l^*{>0}$ a diagonal matrix of load conductances, and $V_l={\rm col} (V_{n_s+1},\dots,V_{n_s+n_l})$, and
\item[(iii)] constant power: $I_l (V_l)= [V_l]^{-1} P_l^*$, with $P_l^*\in \mathds{R}_{<0}^{n_l}$.
\end{itemize}
To refer to 
{the three load cases} above, we will use the indices ``I", ``Z" and ``P" respectively. {The analysis of this paper will focus on the more general case of a parallel combination of the three loads, thus on the case of ``ZIP" loads. {Moreover, additional and stronger statements} results  on the ``ZI" case will be reported.}

Bearing in mind \eqref{resistive.network}, \eqref{consensus.power.control}, and vectorizing the expressions to avoid cluttered formulas, the closed-loop system is 
\be\label{power.consensus}
\begin{bmatrix}
C_s \dot V_s \\ - I_l(V_l) 
\end{bmatrix}
=
-
\begin{bmatrix}
 [V_s]L_c  C_s^{-1} P_s \\ B_l \Gamma B^T V
\end{bmatrix},
\ee
where 
$V={\rm col} (V_s, V_l)$, $C_s ={\rm diag}(C_{1},\dots,C_{n_s})$, 
{$P_s={\rm col} (P_{1},\dots,P_{n_s})$ given by
\begin{equation}
P_s=  [V_s] I_s=  [V_s] (Y_{ss}V_s+Y_{sl}V_l)
\label{eq:source-powers}
\end{equation}
are source power injections and 
\begin{equation}
I_l{(V_l)} = I_l^*-Y_l^* V_l + [V_l]^{-1} P_l^*
\label{eq:load-currents}
\end{equation}
{are  the load currents. }
The interconnected closed-loop DC microgrid is then entirely described by equations \eqref{power.consensus}, \eqref{eq:source-powers}, \eqref{eq:load-currents}.
{An example of a simple closed-loop DC microgrid with two sources and one  constant impedance load is given  in 
Figure \ref{circuit.example}. 
}

%

\begin{remark}({\bf Nonlinear consensus algorithms)}
To compare the algorithm \eqref{consensus.power.control} with {related nonlinear consensus algorithms proposed in the literature, we} neglect the algebraic constraints and the differentiation between sources and loads. This allows us to rewrite \eqref{consensus.power.control} as 
\[
C \dot V = - [V] L_c C^{-1} [V] B\Gamma B^T V.
\]
The weighted power mean consensus algorithms of \cite{bauso.scl06,cortes.aut2008}, on the other hand,  can be written as 
$[W] \dot V = [V]^{1-r} B\Gamma B^T V$, 
where $W$ is vector of weights satisfying $\mathbb{1}^T W=0$ and $r\in \mathbb{R}$. In the special case $r=0$, we get 
\[
[W] \dot V = [V]B\Gamma B^T V,
\]
which is known to converge to the consensus value $V_1^{w_1}\ldots V_n^{w_n}$. The analysis is based on the Lyapunov function
$\sum_{i=1}^n w_i V_i- \prod_{i=1}^n V_i^{w_i}$. 
\\
The nonlinear power consensus algorithm is different in that it uses another  layer of averaging  in addition to the averaging induced by the physical network. This, and the algebraic constraints, requires a different analysis based on physically inspired Lyapunov functions. 
\end{remark}


\begin{figure*}[ht]
  \begin{center}
    \begin{circuitikz}
      \draw (0,0)
      to[I,i=$u_1$] (0,2) 
      to[short] (2,2)
      to[C=$\mathcal{C}_{1}(V_{1})$] (2,0) 
      to[short] (0,0);
      \draw (2,2)
      to[R=$r_1$, i=$I_{1}$] (5,2)
      to[R=$r_3$, i=$I_{3}$] (5,0)
      to[short] (2,0);
      \draw (5,2)
      to[R=$r_2$, i=$-I_{2}$] (8.5,2)
      to[C=$\mathcal{C}_{2}(V_{2})$] (8.5,0) 
      to[short] (5,0);
      \draw (12,0)
      to[I,i=$u_2$] (12,2) 
      to[short] (8,2);
      \draw (8,0)
      to[short] (12,0);
    \end{circuitikz}
    \caption{Circuit considered in Example \ref{exmp2.zip}.\label{circuit.example}}
  \end{center}
\end{figure*}
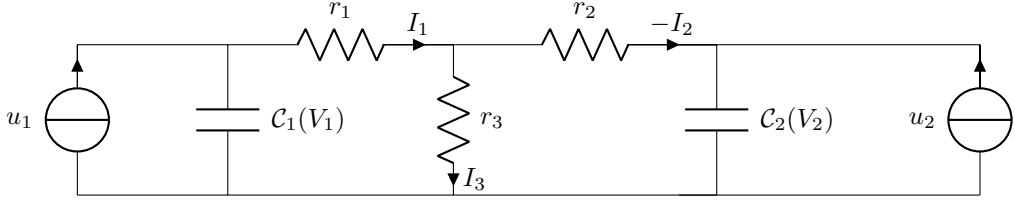

\section{Power consensus algorithm with {ZIP} loads}\label{sec.power.consensus.constant.current.load}
In this section we 
{analyze the closed-loop system \eqref{power.consensus}, \eqref{eq:source-powers}, \eqref{eq:load-currents}.}
We start by studying its equilibria, namely the set of points $V\in \mathds{R}^n_{>0}$ {that} satisfy {\eqref{eq:source-powers}, \eqref{eq:load-currents}, and}
\be\label{power.consensus.constant.current.equilibria}
\begin{bmatrix}
 \mathbb{0} \\ 
 - I_l{(V_l)}
\end{bmatrix}
=
-
\begin{bmatrix}
[V_s]L_c C_s^{-1} P_s \\ B_l \Gamma B^T V
\,.
\end{bmatrix}
\ee


\subsection{{Steady-state characterization}}

{In the following, we show that the equilibria are fully characterized by power balance equations at the sources and current balance equations at the loads, respectively.}

\begin{lem}\label{lem.equilibria}
{\bf (System equilibria)}
The {equilibria of the system  \eqref{power.consensus}, \eqref{eq:source-powers}, \eqref{eq:load-currents} are equivalently characterized by}
\[\ba{l}
\mathcal{E}_{{ZIP}} = \{ V\in \mathds{R}^n_{>0}: {\mathcal{I}_{ZIP}(V)=\mathbb{0}}
\,,\,\mathcal{P}_{{ZIP}}({V})=\mathbb{0} \},
\ea\]
{where $\mathcal{I}_{ZIP}(V)=\mathbb 0$ is the current balance at the loads
\begin{equation*}
\mathcal{I}_{ZIP}(V) = I_l(V_l) - Y_{ll}V_l - Y_{ls} V_s\,,
\end{equation*}
 $\mathcal{P}_{{ZIP}}({V})=0$} depicts the power balance {at} the sources
\[
\mathcal{P}_{{ZIP}}({V})= \underbrace{[ V_s] {Y_{red}}   V_s}_{\ba{c}\text{\tiny network}\\[-3mm] \text{\tiny dissipation} \ea} \!\!+\;    \underbrace{[V_s] Y_{sl} Y_{ll}^{-1} {I_l(V_l)}}_{\ba{c}\text{\tiny load}\\[-3mm] \text{\tiny demands} \ea} \;-\!\!\!\underbrace{P_{s}}_{\ba{c}\text{\tiny source}\\[-3mm] \text{\tiny injections} \ea}\!\!\!\!,
\]
{$Y_{red}$ is the Kron-reduced conductance matrix, ${Y_{ll}^{-1}  Y_{sl}} {I_l(V_l)}$ is the mapping of the {ZIP loads $I_l(V_l)$} to the source buses in the Kron-reduced network as in \eqref{eq:Kron-reduction},} and $P_{s}$ is vector of power injections by the sources written for  $V\in \mathcal{E}_{{ZIP}}$ as
\be\label{Ps}
P_s = - C_s\mathbb{1} \displaystyle\frac{\mathbb{1}^T {I_l(V_l)}}{\mathbb{1}^T [V_s]^{-1}
C_s\mathbb{1}}=: C_s\mathbb{1} p_s^*. 
\ee
\end{lem}

{Observe that the steady-state injections \eqref{Ps} achieve indeed power sharing, and the asymptotic power value $p_s^*$} to which the source power injections converge (in a proportional fashion according to the coefficients $C_i$, $i\in \mathcal{V}_s$) is the total current demand divided by the weighted  sum of the steady-state source voltages. The latter values and those of the load voltages are interestingly entangled by the {power balance at the sources
$\mathcal{P}_{{ZIP}}({V})=\mathbb{0}$ and the current balance equations at the loads $\mathcal{I}_{ZIP}(V) = \mathbb 0$. }

\textbf{Proof.} 
Let $V$ be an equilibrium of {\eqref{power.consensus}, \eqref{eq:source-powers}, \eqref{eq:load-currents}}, that is let $V\in \mathds{R}^n_{>0}$ satisfy  \eqref{power.consensus.constant.current.equilibria}. From the first equation, $ \mathbb{0}=  [V_s]L_c C_s^{-1} P_s$, it immediately follows that $P_s= C_s \mathbb{1}_{{n_s}} p_s^*$ for some scalar $p_s^*$. We rewrite the current balances as
\begin{equation}
\label{power.consensus.current.balance}
\begin{bmatrix}
[V_s]^{-1}C_s \mathbb{1}_{n_s} p_s^*\\ {I_l(V_l)}
\end{bmatrix}
=
\begin{bmatrix}
B_s \Gamma B^T V \\ B_l \Gamma B^T V
\end{bmatrix}.
\end{equation}
Next, we left-multiply \eqref{power.consensus.current.balance} by $[\,\mathbb{1}_{n_s}^T \;\; \mathbb{1}_{n_l}^T]$ to obtain
\[
\mathbb{1}_{n_s}^T [V_s]^{-1} C_s\mathbb{1}_{n_s} p_s^*+\mathbb{1}_{n_l}^T {I_l(V_l)}=0.
\]
{The latter equation can be solved for $p_s^*$ as in \eqref{Ps}.}
From ${I_l(V_l)} = B_l \Gamma B^T V$, 
{we obtain (see \eqref{eq:Kron-reduction}) $\mathcal{I}_{ZIP}(V)=\mathbb 0$ or}
\be\label{algebraic.explicit}
V_l = 
- Y_{ll}^{-1} Y_{ls} V_s+ Y_{ll}^{-1}  {I_l(V_l)}, 
\ee
which replaced in {the first equation of \eqref{power.consensus.current.balance}} returns
\[
Y_{ss} V_s + Y_{sl} (- Y_{ll}^{-1} Y_{ls} V_s+ Y_{ll}^{-1}  {I_l(V_l)}) = [V_s]^{-1}C_s \mathbb{1}_{n_s} p_s^*.
\]
{By rearranging the terms, we arrive at}
\[
{Y_{red}}  V_s + Y_{sl}Y_{ll}^{-1}  {I_l(V_l)}-[V_s]^{-1}C_s \mathbb{1}_{n_s} p_s^*=\mathbb{0},
\] 
which {can be reformulated as $\mathcal{P}_{{ZIP}}({V})=0$ after left-multiplying by $[V_s]$ and bearing in mind \eqref{Ps}.}  The latter and 
\eqref{algebraic.explicit} show that $V\in \mathcal{E}_{{ZIP}}$. 

Conversely, let $V\in \mathcal{E}_{{ZIP}}$. Then the equation $I_l{(V_l)} = B_l \Gamma B^T V$ in \eqref{power.consensus.constant.current.equilibria} is trivially satisfied.   From $\mathcal{P}_{{ZIP}}({V})=0$, and $I_l({V_l}) = B_l \Gamma B^T V$ written as \eqref{algebraic.explicit},  and going backwards through the passages above, we arrive at 
\[
Y_{ss} V_s + Y_{sl} V_l = [V_s]^{-1} C_s\mathbb{1}_{n_s} p_s^*, 
\]
or equivalently at $[V_s] B_s\Gamma B^T V = C_s \mathbb{1}_{n_s} p_s^*$. Hence, the power vector $P_s=[V_s] B_s\Gamma B^T V$ satisfies $L_c C_s^{-1} P_s=\mathbb{0}$, that is, the first equation in  \eqref{power.consensus.constant.current.equilibria}. {Hence, $V\in \mathcal{E}_{ZIP}$ implies that the equilibrium equations \eqref{power.consensus.constant.current.equilibria} are met.}
\qedp

We make the standing assumption that equilibria exist:
\begin{assum}\label{standing}
$\mathcal{E}_{{ZIP}}
\ne \emptyset$.
\end{assum}

\begin{remark}
{\bf (Existence of the equilibria  $\mathcal{E}_{{ZIP}}$)} The analytical investigation of the existence of the equilibria $\mathcal{E}_{{ZIP}}$ is deferred to a future research. This is a topic of interest on its own and similar problems have been dealt with in recent work about the solvability of reactive power flow equations \cite{bolognani.tps16,barabanov,john.resistive.nets,JWSP-FD-FB:14c}. For instance, the problem in \cite{JWSP-FD-FB:14c} boils down to the solution of quadratic algebraic equations of the form $[V_l] Y_{ll} V_l - [V_l] Y_{ll} V_l^*+ Q_l=0$, where $Q_l$ is the vector of constant power load demands and $V_l^*$ is the so called vector of open circuit voltages (again constant). Although  similarities between these equations and the equations $\mathcal{P}_{{ZIP}}(V_s)=\mathbb{0} {=[ V_s] {Y_{red}}   V_s+[V_s] Y_{sl} Y_{ll}^{-1} {I_l(V_l)}+P_{s}}$ could be useful to investigate the nature of the set $\mathcal{E}_{ZIP}$, the non-quadratic nature of  $\mathcal{P}_{ZIP}(V_s)=\mathbb{0}$, {as well as the presence of the additional equations $Y_{ll}^{-1} {I_l(V_l)-V_l} = Y_{ll}^{-1}  Y_{ls} V_s$ pose} additional challenges. {Extra} insights could come from the convex relaxation of the DC power flow equations in the context of optimal {DC power flow} dispatch \cite{lavei.rantzer}. 
\end{remark}

\begin{remark}{\bf (Equilibrium power balance and voltage inequalities)}
To gain further insights into the equilibrium set $\mathcal{E}_{ZIP}$, recall that the vector of power injections is $P ={\rm col} (P_{1},P_{l}) = [V] B\Gamma B^T V$, where $P_{l} = [V_{l}]I_l(V_l)$. Thus, we have the inherent power balance
\begin{equation}
\mathds{1}^T P_s + \mathds{1}^T P_l = V^T B\Gamma B^T V\ge 0
\label{avp}
\end{equation}
implying that the amount of supplied power has to make up for load demands and resistive losses. In the special case of constant power loads, $I_l (V_l)= [V_l]^{-1} P_l^*$, we obtain the  total (or average) power inequality $\mathds{1}^T P_s + \mathds{1}^T P_l^*\geq 0$. Equivalently, after using \eqref{Ps}, we arrive at
\[
- \mathbb{1}^T C_s \mathbb{1} \displaystyle\frac{\mathbb{1}^T [V_l]^{-1} P_l^*}{\mathbb{1}^T [V_s]^{-1} C_s \mathbb{1}} +\mathds{1}^T P_l^*\ge 0\,.
\]  
This inequality can be reformulated as 
\be\label{avi}
\displaystyle\sum_{i\in \mathcal{V}_l} \frac{a_{i}}{V_{i}} \geq \displaystyle\sum_{i\in \mathcal{V}_s} \frac{b_{i}}{V_{i}},
\ee
with $a_i =  P_{l,i}^*/\sum_{i\in \mathcal{V}_l} P_{l,i}^*$ and 
$b_i =  C_{i}/\sum_{i\in \mathcal{V}_s} C_{i}$,
which relates a convex combination of the {reciprocals}  of the voltages at the loads, with a convex combination of the {reciprocals}   of the voltages at the sources, and represents another relation between $V_s, V_l$ in addition to those in \eqref{avp}. The average voltage inequality \eqref{avi} implies that the {reciprocal of the harmonic} average source voltage must be larger than the {reciprocal of the harmonic}  average load voltage so that power can flow from sources to loads.
\end{remark}

{In some special cases reviewed in the  examples below, an explicit characterization of the equilibria can be given.}

\begin{example}\label{exmp2.zip}
Consider the case of two sources ($n_s=2$) and one load ($n_l=1$) as in Figure \ref{circuit.example}, {in which the constant impedance load is replaced by a ZIP load}. 
The equations $\mathcal{P}_{ZIP}(V)=\mathbb{0}$, are in this case
\[\ba{l}
\displaystyle\frac{ \gamma_1 \gamma_2}{\gamma_1+\gamma_2} V_{1}  (V_{1} - V_{2}) -  \frac{ \gamma_1}{\gamma_1+\gamma_2} V_{1}{I_l(V_l)}+\\
\displaystyle \hspace{5cm} {I_l(V_l)}\frac{V_{1} V_{2}}{{V_{1}}+V_{2}} =0\\[4mm]
\displaystyle\frac{ \gamma_1 \gamma_2}{\gamma_1+\gamma_2} V_{1}  (V_{2} - V_{1}) -  \frac{ \gamma_2}{\gamma_1+\gamma_2} V_{2}{I_l(V_l)}+ \\
\displaystyle \hspace{5cm} {I_l(V_l)}\frac{V_{1} V_{2}}{V_{1}+V_{2}}  =0.
\ea\]
We study solutions to the algebraic equations on the curve $V_{1} V_{2} =:c$. The reason for this choice will become clear in Subsection \ref{sec.convergence}.
On such a curve, the equations simplify as
\be\label{equilibria.T.configuration.zip}
\ba{rr}
V_{{1}}^4  - r_2 {I_l(V_l)} V_{{1}}^3 + c r_1 {I_l(V_l)} V_{{1}}-c^2 &=0\\ 
V_{{2}}^4  - r_1 {I_l(V_l)} V_{{2}}^3 + c r_2 {I_l(V_l)} V_{{2}}-c^2 &=0,
\ea
\ee
where $r_i=\gamma_i^{-1}$, $i=1,2$ (the resistance of the transmission line $i$ connecting the source $i$ to the load). \\
{We want to study the solutions of these equations as functions of $I_l(V_l)$. Then} these {can be regarded as} two independent quartic functions for which an analytic, although involved,  expressions of the solutions exist according to the Ferrari-Cardano's formula. 
These expressions simplify if one takes $r_1=r_2$. Then there is a unique positive solution given by $V_1=V_2=\sqrt{c}$, independent of $I_l(V_l)$. The value of $V_l$ is obtained from the algebraic equation 
$0= V_l B_l\Gamma B^T V {-} I_l(V_l)$, solving
\be\label{example.load.voltage}\ba{lcl}
0&=&V_l (-\gamma_1 V_1 - \gamma_2 V_2 +(\gamma_1+\gamma_2) V_l)-  I_l(V_l)\\
&=& 2 \gamma V_l^2 - 2\gamma \sqrt{c} V_l - I_l^*+Y_l^* V_l - [V_l]^{-1} P_l^*\\
&=&2 \gamma V_l^3+ (Y_l^*- 2\gamma \sqrt{c}) V_l^2- I_l^* V_l - P_l^*.
\ea\ee 
In the absence of loads, we have three real roots: a double root at $V_{l}=0$ and a single root at $V_{l}=\sqrt{c} = V_{1}=V_{2}$. Since the roots of a polynomial are continuous in the {parameters}, at most the double root at $V_{l}=0$ can turn to a complex-conjugate root for small loading, and the real root near $\sqrt{c}$ persists.
\end{example}

\subsection{A Lyapunov function {and hidden gradient form}}

We pursue a Lyapunov-based analysis of the stability {of the closed-loop system \eqref{power.consensus}, \eqref{eq:source-powers}, \eqref{eq:load-currents}.} 
Inspired by the Lyapunov analysis of {the} reactive power consensus algorithm in \cite{claudio.nima.ecc16}, we consider the total power dissipated through the network resistors, {$\frac{1}{2} V^T B\Gamma B^T V$,} as the first natural Lyapunov candidate for our analysis, {to which we add  the {power} dissipated through the impedance loads, to obtain {the power losses at passive devices as}
\be\label{J}
J(V)= \frac{1}{2} V^T \left( B\Gamma B^T+\begin{bmatrix} \mathbb{0} & \mathbb{0} \\ \mathbb{0} & {Y_l^*} \end{bmatrix} \right) V. 
\ee
}
Let $\overline V\in \mathcal{E}_{{ZIP}}$, and define $\overline P_s={[\overline V_s] B_s \Gamma B^T \overline V}$ the source power injection corresponding to the equilibrium source voltage ${\overline V}$ (see \eqref{Ps}). 
To cope with the asymmetry in the dynamics of the sources and loads we add to ${J}$ the {terms} 
\[
H(V)= - \overline P_s^T \boldsymbol{\ln} (V_s),
\]
{and} 
\[
{K(V)= - {P_l^*}^T \boldsymbol{\ln}(V_l),}
\]
\noindent which is the way classical power systems transient stability analysis 
absorbs constant power injections \cite{H-DC:11} into a so-called {{\em energy function} defined here as} 
{
\be\label{M}\ba{l}
M(V):=J(V)+H(V)+ K(V)\\
= \frac{1}{2} V^T (B\Gamma B^T+\begin{bmatrix} 0 & 0 \\ 0 & {Y_l^*} \end{bmatrix}) V- \overline P_s^T \boldsymbol{\ln} (V_s) - {P_l^*}^T \boldsymbol{\ln}(V_l).
\ea\ee  
}{The natural ``energy function" \eqref{M} has its minimum at the trivial zero voltage level. To center the function $M$ with respect to a non-trivial equilibrium $\overline V \in \mathcal{E}_{ZIP}$, we use the following Bregman function \cite{claudio.nima.ecc16}}
%
\be\label{calM}
\mathcal{M}(V)= M(V)-M(\overline V)- 
\left.\displaystyle\frac{\partial M}{\partial V}\right|_{V=\overline V}^T(V-\overline V).
\ee
The next result shows a {(perhaps surprising)} gradient relation between the dynamics of {system \eqref{power.consensus}, \eqref{eq:source-powers}, \eqref{eq:load-currents}} and the Bregman function {\eqref{calM}} above:

\begin{lem}\label{gradient.form.constant.current}
{\bf (Gradient dynamics)}
The 
following holds
\be\label{gradient.dynamics.zip}
\begin{bmatrix}
L_c C_s^{-1} P_s \\ B_l \Gamma B^T V - I_l({V_l})
\end{bmatrix}
=
\begin{bmatrix}
L_c [V_s] C_s^{-1} & \mathbb{0}\\
\mathbb{0} & \mathbb{I}_{n_l}
\end{bmatrix}
\displaystyle\frac{\partial {\mathcal{M}}(V)}{\partial V}
\ee
{for all $V\in \mathbb{R}^n_{>0}$. }
Hence the {system \eqref{power.consensus}, \eqref{eq:source-powers}, \eqref{eq:load-currents} can be rewritten as a weighted gradient flow}
\be\label{gradient.dynamics.zip.2}
\ba{rcl}
\begin{bmatrix}
C_s \dot V_s \\ \mathbb{0} 
\end{bmatrix}
&=&
-
\begin{bmatrix}
[V_s]L_c [V_s] C_s^{-1} & \mathbb{0}\\
\mathbb{0} & \mathbb{I}_{n_l}
\end{bmatrix}
\displaystyle\frac{\partial {\mathcal{M}}(V)}{\partial V}.
\ea\ee
\end{lem}

\textbf{Proof.} 
The gradient of the 
function $M(V)$ writes as 
\[
\displaystyle\frac{\partial M}{\partial V} = 
B\Gamma B^T V 
+\begin{bmatrix}
\mathbb{0} \\ Y_l^* V_l
\end{bmatrix}
- \begin{bmatrix} 
[V_s]^{-1}\overline P_s \\ \mathbb{0}
\end{bmatrix}
-
\begin{bmatrix}
\mathbb{0}\\
[V_l]^{-1}{P_l^*}
\end{bmatrix}.
\]
Hence, the  {Bregman  function \eqref{calM}} satisfies
\[\ba{l}
\displaystyle\frac{\partial \mathcal{M}}{\partial V}=\displaystyle\frac{\partial M}{\partial V}-\left.\displaystyle\frac{\partial M}{\partial V}\right|_{{V=\overline V}} = 
\begin{bmatrix} 
B_s \Gamma B^T (V-\overline V)\\
B_l \Gamma B^T (V-\overline V)\\
\end{bmatrix}
+
\\
\begin{bmatrix}
\mathbb{0} \\ Y_l^* (V_l-\overline V_l)
\end{bmatrix}
- \begin{bmatrix} 
([V_s]^{-1}-[\overline V_s]^{-1})\overline P_s \\ 
([V_l]^{-1}-[\overline V_l]^{-1}){P_l^*}
\end{bmatrix}.
%
\ea\]
Bearing in mind the equilibrium condition at the loads
\[
B_l \Gamma B^T \overline V = {I_l{({\overline V}_l)}=} I_l^* - Y_l^* \overline V_l +[\overline V_l]^{-1} P_l^*,
\]
and replacing it in the second line of the identity {above describing $\displaystyle{\partial \mathcal{M}}/{\partial V}$},  
we obtain 
\[\ba{rcl}
\displaystyle\frac{\partial \mathcal{M}}{\partial V_l} &=&
B_l \Gamma B^T V + Y_l^*  V_l 
-[V_l]^{-1} {P_l^*} - I_l^*
\\
&=&
B_l \Gamma B^T V 
-I_l(V_l),
\ea\]
which {equals precisely} the second equation in \eqref{gradient.dynamics.zip}. 

Analogously, for {the first line ${\partial \mathcal{M}}/{\partial V_s}$,} we write
\be\label{dM.dVs}
\ba{rcl}
\displaystyle\frac{\partial \mathcal{M}}{\partial V_s} &=& B_s \Gamma B^T (V-\overline V)-([V_s]^{-1}-[\overline V_s]^{-1})\overline P_s 
\\
 &=& [V_s]^{-1} P_s -[\overline V_s]^{-1} \overline P_s -([V_s]^{-1}-[\overline V_s]^{-1})\overline P_s \\
 &=& [V_s]^{-1} (P_s - \overline P_s),
\ea
\ee
where to write the second equality we have used the identities $P_s = [V_s] B_s \Gamma B^T V$ and $\overline P_s = [\overline V_s] B_s \Gamma B^T \overline V$. 

Now note that 
\[
P_s= [V_s] \displaystyle\frac{\partial \mathcal{M}}{\partial V_s} +\overline P_s
\]
and, multiplying both sides by $L_c C_s^{-1}$, we obtain
\[\ba{rcl}
L_c C_s^{-1} P_s &=& L_c C_s^{-1} [V_s] \displaystyle\frac{\partial \mathcal{M}}{\partial V_s} +L_c C_s^{-1} \overline P_s\\[3mm]
&=& L_c C_s^{-1} [V_s] \displaystyle\frac{\partial \mathcal{M}}{\partial V_s}, 
\ea\]
having exploited that $\overline V \in \mathcal{E}_{ZIP}$ implies  $\overline P_s= C_s \mathbb{1} p_s^*$. The identity 
$L_c C_s^{-1} P_s = L_c C_s^{-1} [V_s] \displaystyle\frac{\partial \mathcal{M}}{\partial V_s}$ is the first equation in 
\eqref{gradient.dynamics.zip}.

In view of the dynamics {\eqref{power.consensus}, \eqref{eq:source-powers}, \eqref{eq:load-currents}}, one immediately realizes that 
\[
\begin{bmatrix}
L_c C_s^{-1} P_s \\ B_l \Gamma B^T V - I_l(V_l)
\end{bmatrix}
= 
\begin{bmatrix}
- [V_s]^{-1} C_s \dot V_s \\ \mathbb{0}
\end{bmatrix},
\]
{showing the identity \eqref{gradient.dynamics.zip.2} which concludes the proof.}
\qedp

%

\subsection{Convergence of solutions}\label{sec.convergence}

The particular form of the dynamics {\eqref{power.consensus}, \eqref{eq:source-powers}, \eqref{eq:load-currents}}  elucidated in Lemma \ref{gradient.form.constant.current} permits a straightforward analysis of the convergence properties of the solutions. 

\begin{thm}\label{thm.constant.current}
{\bf (Main result)}
{Assume that there exists $\overline V\in \mathcal{E}_{{ZIP}}$ such that 
\be\label{main.cond}
Y_{ll}+ Y_l^* +[\overline V_l]^{-2} [P_l^*]- Y_{ls} (Y_{ss}+[\overline V_s]^{-2} [\overline P_s])^{-1} Y_{sl}>0
\ee
Then there exists a compact sublevel set $\Lambda_{ZIP}$ of   $\mathcal{M}$ contained in $\mathds{R}^n_{>0}$ such that 
any solution 
}
to {\eqref{power.consensus}, \eqref{eq:source-powers}, \eqref{eq:load-currents}} that originates from initial conditions $V(0)$ belonging to  $\Lambda_{{ZIP}}$ exists,  always remain in $\Lambda_{ZIP}$ {with strictly positive voltages for all times,} 
and asymptotically converges to the set 
$\mathcal{E}_{{ZIP}}\cap \Lambda_{{ZIP}}\cap \mathcal{V}_{{ZIP}}$,  
where {$\mathcal{V}_{ZIP}$ specifies the preserved weighted geometric mean of the source voltages
\[\ba{ll}
 \mathcal{V}_{ZIP}:= &\{(V_s,V_l)\in \Lambda_{ZIP}:  \mathcal{I}_{ZIP}(V) = \mathbb 0\,,\,\\
 & V_1^{C_{1}}\cdot\ldots \cdot V_{n_s}^{C_{n_s}}=V_1^{C_{1}}(0)\cdot\ldots \cdot V_{n_s}^{C_{n_s}}(0)
 \}.
\ea\]}
\end{thm}

\begin{remark}{\bf (Interpretation of the main condition)}
The main condition \eqref{main.cond} guarantees regularity of the algebraic equations and stability of the solutions. Its role is revealed when converting the constant power loads and the asymptotically constant power injections at the sources to the equivalent impedances $[\overline V_l]^{-2} [P_l^*]$ and $[\overline V_s]^{-2} [\overline P_s]$. In this case, the equivalent conductance matrix in the steady-state current-balance equations \eqref{resistive.network} read as
\begin{equation}
\label{eq:Yeq}
	Y_{eq} = 
\begin{bmatrix}
Y_{ss} & Y_{sl}\\
Y_{ls} & Y_{ll}
\end{bmatrix}
+ \begin{bmatrix} 
[V_s]^{-2} [\overline P_s] &  \mathbb{0} \\ 
 \mathbb{0} & [V_l]^{-2}[P_l^*] + Y_{l}^{*}
\end{bmatrix}
\,.
\end{equation}
By a Schur complement argument, observe that $Y_{eq}$ is a well-defined (i.e., positive definite) conductance matrix if and only if the main condition \eqref{main.cond} holds.
\end{remark}


\textbf{Proof.} 
{\it Existence and boundedness of solutions.} Observe {first} that 
\be\label{calS.hessian}
\frac{\partial^2 \mathcal{{M}}}{\partial V^2} = B\Gamma B^T 
{
+\begin{bmatrix}
\mathbb{0}  & \mathbb{0} \\ 
\mathbb{0}  & Y_l^* 
\end{bmatrix}
+ \begin{bmatrix} 
[V_s]^{-2} [\overline P_s] &  \mathbb{0} \\ 
 \mathbb{0} & [V_l]^{-2}[P_l^*]
\end{bmatrix},
 }
%
\ee
Let $\overline V {>0}$ be an equilibrium of the system, i.e.,~$\overline V\in \mathcal{E}_{{ZIP}}$.
Since $I_l^*, P_l^*\in\mathds{R}^{n_l}_{<0}$, {and $\overline V>0$,} the {steady-state} power injection at the sources satisfies $\overline P_s\in \mathds{R}^{n_s}_{>0}$ by \eqref{Ps}. {Hence, $[V_s]^{-2} [\overline P_s]>0$  is positive definite.}
{Then} 
the  Bregman function ${\mathcal{M}}$ has an isolated minimum at the equilibrium $\overline V$, {in view of \eqref{main.cond}, \eqref{calS.hessian} and a standard Schur complement argument.} Then there exists a compact sublevel set $\Lambda_{{ZIP}}$ of $\mathcal{{M}}$ around the equilibrium $\overline V$  contained in the positive orthant. {Without loss of generality this compact sublevel set can be taken 
so that all the solutions to  {\eqref{power.consensus}, \eqref{eq:source-powers}, \eqref{eq:load-currents}}  that originate here locally exist.} 

{The algebraic equations \eqref{eq:load-currents} written as in Lemma~\ref{lem.equilibria} are
\begin{align*}
\mathbb 0 =& \mathcal{I}_{ZIP}(V) = I_l(V_l) - Y_{ll}V_l - Y_{ls} V_s  
\\
=& I_l^* {-Y_l^* V_l + [V_l]^{-1} P_l^*} - Y_{ls} V_s- Y_{ll} V_l
\,.
\end{align*}
To study local solvability of these equations, we analyze 
\begin{equation*}
\label{reg.cond}
\frac{\partial \mathcal{I}_{ZIP}}{\partial V_l}= - \left( Y_{ll} + Y_l^*+ [V_l]^{-2} [P_l^*] \right) \,.
\end{equation*}}%
 {In view of 
 \eqref{main.cond}}, nonsingularity of {${\partial \mathcal{I}_{ZIP}}/{\partial V_l}$} and therefore regularity of the algebraic condition holds {in a neighborhood of {$\overline V \in \Lambda_{ZIP}$ from the implicit function theorem \cite{RA-JEM-TSR:88}.} The  sublevel  set $\Lambda_{ZIP}$ can be taken sufficiently small such that it is  contained in the neighborhood of regularity for the algebraic equations, thus showing the claim that solutions starting from $\Lambda_{ZIP}$ locally exist in time, see \cite[Theorem 1]{DJH-IMYM:90} and \cite[Lemma 2.3]{Johannes.florian.ecc16}.
 }

When computed along these solutions, ${\mathcal{M}(V(t))}$ satisfies  
\[
\dot{\mathcal{M}}(V(t))= 
\left.\frac{\partial \mathcal{{M}}}{\partial V_s}\right|_{V=V(t)}^T \dot V_s(t) +\left.\frac{\partial \mathcal{{M}}}{\partial V_l}\right|_{V=V(t)}^T \dot V_l(t).
\]
Notice 
that, {by the algebraic constraint \eqref{gradient.dynamics.zip},} 
\[
\left.\frac{\partial \mathcal{{M}}}{\partial V_l}\right|_{V=V(t)}  = 
B_l\Gamma B^T {V(t)}-I_l(V_l{(t)})=\mathbb{0}
\]
{
for all $t$ for which a solution exists. Hence, we arrive at}
\begin{multline*}
\dot{\mathcal{M}}(V(t)) =
\displaystyle \left.\frac{\partial \mathcal{{M}}}{\partial V_s}\right|_{V=V(t)}^T \dot V_s(t) \\
=
- \displaystyle \left.\frac{\partial \mathcal{{M}}}{\partial V_s}\right|_{V=V(t)}^T 
C_s^{-1} 
[V_s] 
L_c 
{[V_s]}
C_s^{-1}\displaystyle 
\left. 
\frac{\partial \mathcal{M}}{\partial V_s}\right|_{V=V(t)} \leq 0,
\end{multline*}
{where the second equality holds because of \eqref{gradient.dynamics.zip.2}.
The inequality above  shows that} $\mathcal{M}(V(t))$ is a {non-increasing} function of time. By the compactness of the sublevel set around $\overline V$,  
the solutions are bounded, exist and belong to $\Lambda_{{ZIP}}$ for all times. {Thus, among others the voltages stay positive for all times.}

{\it Convergence.} 
Exploiting the regularity of the algebraic equation, the DAE system can be reduced to an ODE system and then the standard LaSalle invariance principle for ODE can be used to infer convergence, see also \cite{Johannes.florian.ecc16}. We argue as follows. 
%
%
Any solution $(V_s, V_l)$ to the DAE system {\eqref{power.consensus}, \eqref{eq:source-powers}, \eqref{eq:load-currents}} {originating in $\Lambda_{ZIP}$} is such that its component  $V_s$ is a solution to the  system of ODE 
\be\label{ode}
\dot V_s = -C_s^{-1} [V_s] L_c [V_s] C_s^{-1}
{
(Y_{ss} V_s+Y_{sl} \delta(V_s)),
}
\ee
{where the map $V_l=\delta(V_s)$ denotes the 
solution of the  algebraic equation {$ \mathcal{I}_{ZIP}(V)=\mathbb{0}$} in $\Lambda_{ZIP}$. }
%
%
%
Define 
\be\label{calN}
\mathcal{N}(V_s) := \mathcal{\mathcal{M}}(V_s, \delta(V_s))
\ee
and observe that 
\[\ba{rcl}
\dot{\mathcal{N}}(V_s(t)) &= &
\displaystyle
\left.\frac{\partial \mathcal{{M}}}{\partial V_s}\right|_{{\tiny \ba{l} V_s= V_s(t)\\[-1.5mm] V_l=\delta(V_s(t))\ea}}^T \dot V_s(t) 
+\left.\frac{\partial \mathcal{{M}}}{\partial V_l}\right|_{{\tiny \ba{l} V_s= V_s(t)\\[-1.5mm] V_l=\delta(V_s(t))\ea}}^T
\cdot \\
&& \displaystyle\qquad\qquad \left.\frac{\partial \delta}{\partial V_s}\right|_{V_s= V_s(t)}
\dot V_s(t)\\
&=& \displaystyle\left.\frac{\partial \mathcal{{M}}}{\partial V_s}\right|_{{\tiny \ba{l} V_s= V_s(t)\\[-1.5mm] V_l=\delta(V_s(t))\ea}}^T \dot V_s(t), 
\ea
\]
since 
\[\ba{rcl}
\displaystyle\left.\frac{\partial \mathcal{{M}}}{\partial V_l}\right|_{{\tiny \ba{l} V_s= V_s(t)\\[-1.5mm] V_l=\delta(V_s(t))\ea}}
&=& Y_{ls} V_s(t)+Y_{ll} \delta(V_s(t))-I_l( \delta(V_s(t)))\\
&=& Y_{ls} V_s(t)+Y_{ll} V_l(t)-I_l(V_l(t))
=\mathbb{0}
\ea\]
where the second equality holds because $V_l(t)=\delta(V_s(t))$ on $\Lambda_{ZIP}$ and the third equality because of 
the algebraic equation in {\eqref{power.consensus}, \eqref{eq:source-powers}, \eqref{eq:load-currents}}.
It then follows that
%
%
\be\label{dot.calS}\ba{ll}
\dot{\mathcal{N}}(V_s) &= (P_s-\overline P_s)^T [V_s]^{-1}\dot V_s\\
& = - (P_s-\overline P_s)^T   C_s^{-1}  L_c C_s^{-1} P_s\\
&= -  P_s^T   C_s^{-1}  L_c C_s^{-1} P_s \leq 0, 
\ea\ee
where the first equality descends from \eqref{dM.dVs}, the second from {\eqref{power.consensus}}, and the third from \eqref{Ps}.

Since $V_s$ is bounded, then the standard La Salle invariance principle for {ODEs} yields convergence of $V_s$ to the largest invariant set where $L_c C_s^{-1}  P_s=\mathbb{0}$. {Moreover, since the solutions evolve in $\Lambda_{ZIP}$}, since they satisfy the algebraic equations, and since $L_c  C_s^{-1}  P_s=\mathbb{0}$, we have from Lemma \ref{lem.equilibria} that at steady state $(V_s, V_l)\in \mathcal{E}_{{ZIP}}$.
Since $(V_s, V_l)$ is a solution to {\eqref{power.consensus}, \eqref{eq:source-powers}, \eqref{eq:load-currents}}  that remains in $\Lambda_{{ZIP}}$, convergence to the set $\mathcal{E}_{{ZIP}}\cap \Lambda_{{ZIP}}$ is inferred.  Moreover, the quantity $V_1\cdot\ldots \cdot V_{n_s}$ is conserved, namely $V_1(t)\cdot\ldots \cdot V_{n_s}(t)=V_1(0)\cdot\ldots \cdot V_{n_s}(0)$ for all $t$. In fact,  by \eqref{ode},
$C_s\frac{d}{dt} \boldsymbol{\ln} V_s = - L_c [V_s] C_s^{-1} 
{
(Y_{ss} V_s+Y_{sl} \delta(V_s))
}
$, and therefore $\frac{d}{dt} \mathbb{1}^T C_s \boldsymbol{\ln} V_s=0$.
The thesis then follows. 
\qedp
\begin{example}\label{exmp4}
Consider again the case of two sources ($n_s=2$) and one load ($n_l=1$) connected in a ``T" configuration, as in Example {\ref{exmp2.zip}}.  If $C_s=\mathbb{I}_2 {C}$, {for some positive real number $C$,} the result above shows that on the convergence set $\mathcal{E}_{{ZIP}} \cap \Lambda_{{ZIP}}\cap \mathcal{V}_{{ZIP}}$, 
$V_{1} V_{2} = V_{1}(0) V_{2}(0)=:c$ for all $t\ge 0$. Hence, as discussed in  Example {\ref{exmp2.zip}}, 
the expression of the (real and positive) solution to the {equations \eqref{equilibria.T.configuration.zip}} takes on a particularly simple form, namely
$V_{1}= V_{2} =\sqrt{c}=\sqrt{V_{{1}}(0) V_{{2}}(0)}$.
It follows that any point on $\mathcal{E}_{{ZIP}}\cap \mathcal{V}_{{ZIP}}$  is such that each source voltage is the geometric mean of the initial voltage sources. 
Accordingly, {the load voltage  must satisfy \eqref{example.load.voltage}.}
\end{example}

\begin{remark}\label{rem.capacitor.current.loads}
{\bf (Capacitors at the loads)} If loads are interconnected to the network via capacitors, the load equations are modified as
\[
C_l \dot V_l = -{I_l(V_l)} + B_l\Gamma B^T V.
\]
Notice that the equilibria of the system remain the same. 
Bearing in mind {\eqref{gradient.dynamics.zip},} the load dynamics {read as}  
\[
C_l \dot V_l =  -\displaystyle\frac{\partial \mathcal{{M}}}{\partial V_l}. 
\]
It follows that 
\[
\dot{\mathcal{{M}}}= - \frac{\partial \mathcal{{M}}}{\partial V_s}^T  
C_s^{-1} 
[V_s] 
L_c 
{[V_s]}
C_s^{-1}\displaystyle \frac{\partial \mathcal{{M}}}{\partial V_s}- \frac{\partial \mathcal{{M}}}{\partial V_l}^T C_l^{-1} \frac{\partial \mathcal{{M}}}{\partial V_l},
\]  
and one can infer convergence to the set $\mathcal{E}_{ZIP}\cap \Lambda_{ZIP}\cap \mathcal{V}_{ZIP}$ similarly as for the differential-algebraic model. 
\end{remark}

\begin{remark}
{\bf (Constant voltage buses)} Similarly as in \cite[Remark 3.3]{zhao.dorfler.aut15}, one can consider also consider voltage-controlled buses. For example, consider the scenario of all load buses having constant (not necessarily identical) voltages $\overline V_l$ (see \cite{zhao.dorfler.aut15} for a discussion on this load condition). More precisely, a controller adjust the current injection $I_l$ depending on $V_s$ to maintain the value of the voltage at the constant level $\overline V_l$ so that system \eqref{power.consensus} {reads} as%
\begin{subequations}\label{power.consensus.constant.voltage}%
\begin{align}
	C_s \dot V_s &= -[V_s]L_c C_s^{-1} [V_s] (Y_{ss} V_s+Y_{sl} \overline V_l)
	\label{power.consensus.constant.voltage-1}\\
- I_l &= Y_{ls} V_s+Y_{ll} \overline V_l. 
\label{power.consensus.constant.voltage-2}
\end{align}%
\end{subequations}%
The only relevant equations for stability of  \eqref{power.consensus.constant.voltage} are the ordinary differential equations \eqref{power.consensus.constant.voltage-1} driven by the constant term $\overline V_l$. We study their stability using a similar  Lyapunov argument as before. 
Since $\overline V_l$ is now constant, we consider a simplified version of the  function $\mathcal M$, namely $\tilde{\mathcal M}(V_{s}) = \tilde{\mathcal J}(V_{s}) + \mathcal H(V_{s})$, where 
$$\mathcal{\tilde J}(V_s)=\frac{1}{2} (V_s-\overline V_s)^T Y_{ss} (V_s-\overline V_s)$$ are the (shifted) network losses so that 
$
	\frac{\partial \mathcal{\tilde{J}}}{\partial V_s} =
	Y_{ss} (V_s-\overline V_s)
	=
(Y_{ss} V_s+Y_{sl} \overline V_l)-(Y_{ss}\overline V_s+Y_{sl}\overline V_l)
=
[V_s]^{-1}P_s - [\overline V_s]^{-1}\overline P_s\,.
$
Together with $\mathcal{H} ({V_s})= - \overline P_s^T \boldsymbol{\ln} (V_s)+  \overline P_s^T \boldsymbol{\ln} (\overline V_s)+ \overline P_s^T [\overline V_s]^{-1} (V_s-\overline V_s)$, 
we obtain that $\frac{\partial\mathcal{{\mathcal{\tilde M}}}}{\partial V_s}=[V_s]^{-1}(P_s-\overline P_s)$ and thus $C_s \dot V_s = -[V_s]L_c C_s^{-1}[V_s]  \frac{\partial\mathcal{\tilde{M}}}{\partial V_s}$. The convergence analysis of the solutions of the system \eqref{power.consensus.constant.voltage} is now analogous to the proof of Theorem \ref{thm.constant.current}.  
%
%
\end{remark}

\subsection{The case of ZI loads}

In the case of ZI loads the {previous results} can be strengthened. First, the set of equilibria can be {characterized} by two systems of equations, one depending on the source voltages only and the other one allowing for a straightforward calculation of the load voltages once the source voltages are determined. Second, the convergence result can be established without any extra {condition on the equivalent conductance matrix in \eqref{eq:Yeq}.}
Finally, the convergence is to a point rather than to a set.

{The first result we present concerns the set of equilibria, which follows by adapting the proof of Lemma~\ref{lem.equilibria}.}

\begin{lem}\label{lem.equilibria.zi.loads}
{\bf (Equilibria for ZI loads)} 
The set of equilibria  of system {\eqref{power.consensus}, \eqref{eq:source-powers}, \eqref{eq:load-currents}} with $I_l(V_l) = I_l^* - Y_l^* V_l$
is
\[\ba{l}
\mathcal{E}_{{ZI}} = \{ V\in \mathds{R}^n_{>0}: \mathcal{P}_{{ZI}}({V_s})=\mathbb{0}, \\ 
\qquad\qquad\qquad
V_l = (Y_{ll}+ Y_l^*)^{-1} (I_l^* - Y_{ls} V_s)\},
\ea\]
where $\mathcal{P}_{{ZI}}(V_s)$ depicts the power balance {at} the sources
\[
\mathcal{P}_{{ZI}}({V})= \underbrace{[ V_s]  {\hat Y_{red}}   V_s}_{\ba{c}\text{\tiny network}\\[-3mm] \text{\tiny dissipation} \ea} \!\!+\;    \underbrace{[V_s] Y_{sl} (Y_{ll}+ Y_l^*)^{-1} {I_l^*}}_{\ba{c}\text{\tiny load}\\[-3mm] \text{\tiny demands} \ea} \;-\!\!\!\underbrace{P_{s}}_{\ba{c}\text{\tiny source}\\[-3mm] \text{\tiny injections} \ea}\!\!\!\!,
\]
 ${\hat Y_{red}}  = Y_{ss}- Y_{sl} (Y_{ll}+ Y_l^*)^{-1} Y_{ls}$ {is the Kron-reduced conductance matrix} that also absorbed the constant impedance loads, 
 and $P_{s}$ is vector of power injections by the sources written for  $V\in \mathcal{E}_{{ZI}}$ as
 $C_s\mathbb{1} p_s^*$, with 
\[
p_s^* = - \mathbb{1}^T\displaystyle\frac{I_l^*-Y_l^* (Y_{ll}+ Y_l^*)^{-1}(I_l^*-Y_{ls} V_s)}
{\mathbb{1}^T [V_s]^{-1}
C_s\mathbb{1}}. 
\]

\end{lem}

We remark that in {the ZI case the} equations $\mathcal{P}_{ZI}(V_s)=\mathbb{0}$ depend on the source voltages only, and once a solution to it is determined, the corresponding voltages at the loads are obtained {as} $V_l = (Y_{ll}+ Y_l^*)^{-1} (I_l^* - Y_{ls} V_s)$ {thereby explicitly solving previous $\mathcal{I}_{ZI}(V) = \mathbb 0$.} 


Our second result concerns the convergence of the dynamics. In the case of ZI loads, convergence can be established without the definiteness condition on the equivalent conductance matrix $Y_{eq}$ in \eqref{eq:Yeq}. Indeed, for $P_{l}^{*}=\mathbb 0$, the condition \eqref{main.cond} is automatically satisfied. Before, this condition was needed to certify strict convexity of the shifted Lyapunov function $\mathcal{M}$ (see \eqref{calS.hessian}) as well as the regularity of the algebraic equation $\mathcal{I}_{ZI}(V) = \mathbb 0$. Additionally,
the {limit} set in case of  ZI loads is $\mathcal{E}_{ZI}\cap{\Lambda}_{ZI}\cap \mathcal{V}_{ZI}$, where the set of equilibria  $\mathcal{E}_{ZI}$ {is} characterized in Lemma \ref{lem.equilibria.zi.loads}, ${\Lambda}_{ZI}$ is a sublevel set associated with the Lyapunov function $\mathcal{M}$ with $P_l^*=\mathbb{0}$, and the set $\mathcal{V}_{ZI}$ is defined as 
 \[
 \ba{ll}
 \mathcal{V}_{{ZI}}:= &\{(V_s,V_l)\in \Lambda_{{ZI}}:\\
 & V_1^{C_{1}}\cdot\ldots \cdot V_{n_s}^{C_{n_s}}=V_1^{C_{1}}(0)\cdot\ldots \cdot V_{n_s}^{C_{n_s}}(0), \; \\
 & {  V_l = (Y_{ll}+Y_l^*)^{-1} (I_l^*-Y_{ls} V_s)}\}.
\ea\]
%
%
%
%
{Finally}, a stronger convergence result {can be established, namely} any trajectory converges to a point depending on the initial condition. This can be formalized as follows:

\begin{thm}\label{cor.constant.current}
{\bf (Point convergence)}
The solutions to {\eqref{power.consensus}, \eqref{eq:source-powers}, \eqref{eq:load-currents}} {with $P_l^*=\mathbb{0}$} which originate from any initial condition $V(0)$ belonging to a {sublevel} set $\Lambda_{{ZI}}$ of {the shifted Lyapunov function $\mathcal{M}$ \eqref{calM} with  $P_l^*=\mathbb{0}$} contained in $\mathds{R}^n_{>0}$ always remain in $\Lambda_{{ZI}}$ and converge to an asymptotically stable equilibrium belonging to 
$\mathcal{E}_{{ZI}}\cap \Lambda_{{ZI}}\cap \mathcal{V}_{{ZI}}$. 
\end{thm}

\textbf{Proof.} 
{First of all we observe that the proof of Theorem \ref{thm.constant.current} holds for the case of ZI loads (it suffices to set $P_l^*=\mathbb{0}$ and $I_l(V_l)= I_l^*- Y_l^* V_l$ throughout the proof).} {As an additional feature of ZI loads (to be exploited below) we can explicitly construct $\delta (V_s) = (Y_{ll}+Y_l^*)^{-1}(I_l^* - Y_{ls} V_s)$.}%

From the proof of Theorem \ref{thm.constant.current} {(specialized to the case of ZI loads),} it is known that any solution $V_s$ of the ODE {\eqref{ode}} is bounded. By Birckhoff's Lemma (\cite[Lemma 3.1]{Khalil}) the positive limit set $\Omega(V_s)$ associated with a solution $V_s(t)$ is non-empty, compact, and invariant. Moreover, it is contained in $\mathcal{E}_{{ZI}}\cap \Lambda_{{ZI}}\cap \mathcal{V}_{{ZI}}$. We would like to prove that  $\Omega(V_s)$ is a singleton. To this end, and similarly to \cite{claudio.nima.ecc16} we appeal to  \cite[Proposition 4.7]{haddad.book}, which states that if the positive limit set $\Omega(V_s)$ {of a trajectory} contains a Lyapunov stable equilibrium $\overline V_s$, then $\Omega(V_s)=\{\overline V_s\}$. To see this first notice that $\overline V_s$ being in $\Omega(V_s)$ and hence in $\mathcal{E}_{{ZI}}\cap \Lambda_{{ZI}}\cap \mathcal{V}_{{ZI}}$, it is indeed an equilibrium of the system. Thus, following {\eqref{calN},} one can construct a shifted  function ${\mathcal{N}}(V_s)$ associated to $\overline V_s$.  The explicit expression of ${\mathcal{N}}(V_s)$ is given by
\begin{multline*}%
{\mathcal{N}}(V_s)
=
- \overline P_{s}^T \boldsymbol{\ln} (V_s)+  \overline P_{s}^T \boldsymbol{\ln} (\overline V_{s})+
\overline P_{s}^T [\overline V_{s}]^{-1} (V_s-\overline V_{s})
\\+
\frac{1}{2} \begin{bmatrix} 
V_s-\overline V_{s}\\
{\delta}(V_s)-{\delta}(\overline V_{s})
\end{bmatrix}
^T 
\begin{bmatrix}
Y_{ss} & Y_{sl}\\
Y_{ls} & Y_{ll} +  Y_l^*
\end{bmatrix}
\begin{bmatrix} 
V_s-\overline V_{s}\\
{\delta}(V_s)-{\delta}(\overline V_{s})
\end{bmatrix}.
\end{multline*}
The gradient of  $\mathcal{{N}}(V_s)$ is given by
\begin{align*}
\displaystyle\frac{\partial \mathcal{{N}}}{\partial V_s}
&= 
- [V_{s}]^{-1} \overline P_{s} + [\overline V_{s}]^{-1} \overline P_{s} +
\\& \left(Y_{ss} + Y_{sl}\frac{\partial  {\delta}}{\partial V_s}\right)^T (V_s-\overline V_{s})
+\\
&\left(Y_{ls} + (Y_{ll} + Y_{l}^{*}) \frac{\partial  {\delta}}{\partial V_s} \right)^T({\delta}(V_s)-{\delta}(\overline V_{s}))\,.
\end{align*}
Since $\frac{\partial  {\delta}}{\partial V_s} = -  (Y_{ll}+Y_l^*)^{-1}Y_{ls}$, the last summand above vanishes. 
With the shorthand ${\hat Y_{red}}  = Y_{ss} - Y_{sl}(Y_{ll}+Y_l^*)^{-1} Y_{ls}$, the gradient simplifies as}
\[\ba{l}
\displaystyle\frac{\partial \mathcal{{N}}}{\partial V_s} = 
- [V_{s}]^{-1} \overline P_{s} + [\overline V_{s}]^{-1} \overline P_{s} + {\hat Y}_{red} (V_s-\overline V_s)
\ea\]
{Note that the gradient $\frac{\partial \mathcal{{N}}}{\partial V_s}$ vanishes if $V_{s} = \overline V_{s}$ and $\mathcal{{N}}$ has a strict local minimum at $\overline V_{s}$ since}
\begin{align*}
\displaystyle\frac{\partial^2 \mathcal{{N}}}{\partial V_s^2}
&=
{\hat Y}_{red} + [V_s]^{-2}\overline P_{s}.
\end{align*}
By \eqref{dot.calS}, $\dot{\mathcal{{N}}}\le 0$, and these two properties (properness and {the} nonnegative time derivative) show that $\overline V_{s}$  is a Lyapunov stable equilibrium. {Therefore,  $\Omega(V_s)=\{\overline V_s\}$, and} the solution $V_s(t)$ converges to an equilibrium point. Because $V_s(t)$ is the $V_s$ component of the solution to the DAE, and since $V_l$ satisfies {$V_l=\delta (V_s) = (Y_{ll}+Y_l^*)^{-1}(I_l^* - Y_{ls} V_s)$} we also see that the solution  $(V_s(t), V_l(t))$ of the DAE {\eqref{power.consensus}, \eqref{eq:source-powers}, \eqref{eq:load-currents}} converges to a point in $\mathcal{E}_{{ZI}}\cap \Lambda_{{ZI}}\cap \mathcal{V}_{{ZI}}$.  
Since this equilibrium point is Lyapunov stable by \eqref{calS.hessian} {(with $P_l^*=\mathbb{0}$)} and {\eqref{dot.calS}, the limit point} is also  asymptotically stable. 
\qedp

\section{Simulations}\label{simulations}

\begin{figure}
\includegraphics[width=0.5\textwidth]{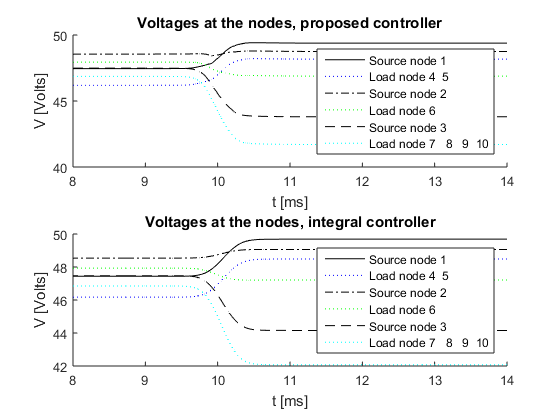}
\caption{Voltage plots of the simulation}
\label{f:simvoltage}
\end{figure}

\begin{figure}
\includegraphics[width=0.5\textwidth]{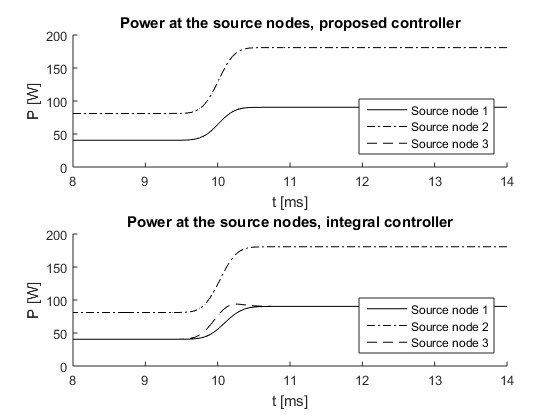}
\caption{Power plots of the simulation}
\label{f:simpower}
\end{figure}

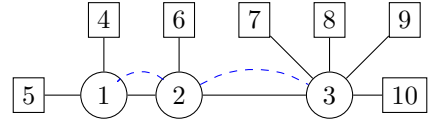
\begin{figure}\centering
\begin{tikzpicture}
	\node [draw, circle] (area1) at (1, 0) {1};
	\node [draw, circle] (area2) at (2, 0) {2};
	\node [draw, circle] (area3) at (4, 0) {3};
	\node [draw, rectangle] (area4) at (1, 1) {4};
	\node [draw, rectangle] (area5) at (0, 0) {5};
	\node [draw, rectangle] (area6) at (2, 1) {6};
	\node [draw, rectangle] (area7) at (3, 1) {7};
	\node [draw, rectangle] (area8) at (4, 1) {8};
	\node [draw, rectangle] (area9) at (5, 1) {9};
	\node [draw, rectangle] (area0) at (5, 0) {10};
	\draw (area1) -- (area2);
	\draw (area2) -- (area3);
	\draw (area1) -- (area4);
	\draw (area1) -- (area5);
	\draw (area2) -- (area6);
	\draw (area3) -- (area7);
	\draw (area3) -- (area8);
	\draw (area3) -- (area9);
	\draw (area3) -- (area0);
	\draw [rounded corners=3ex, blue, dashed] (area1) -- (1.5,0.5) -- (area2);
	\draw [rounded corners=4ex, blue, dashed] (area2) -- (3  ,0.5) -- (area3);
\end{tikzpicture}
\caption{The node network used for the simulations.
Sources are depicted as circles, loads as rectangles.
Solid lines denote the interconnecting lines, while dashed blue lines represent the communication graph used by the controllers.}
\label{f:simnetwork}
\end{figure}

In this section, we present simulation results comparing the proposed control strategy to an averaging-based control method.
We use an example network obtained from \cite{belk2016stability}.
The network topology is sketched in Fig.~\ref{f:simnetwork}, and the physical parameters are given in Table \ref{t:simparams}.
It can be checked that condition \eqref{main.cond} is satisfied for this network. 
As in the reference experiment, there are seven constant power loads, five of which are initially turned off and are turned on gradually between \num{9.5} and \SI{10.5}{\milli\second}. This means that there is a gradual increase of the total power load from \SI{70}{\watt} to \SI{245}{\watt}.  
We simulate both the proposed control strategy 
(5), and the distributed averaging integral controller \eqref{dapi} for a comparison. 
The power measured at the source nodes is shown for both control strategies in Fig.~\ref{f:simpower}. As predicted by the analysis,  at steady state proportional power sharing is achieved by the power sources in conformity with  \eqref{pps}. We also observe that the two controllers perform similarly, only a slight overshoot for the integral controller at the  power source 2 can be observed. The voltage evolution both at the sources and at the loads is depicted in Fig.~\ref{f:simvoltage}. 


\begin{table}\centering
\begin{tabular}{ll}
\hline\hline
Parameter & Value \\
\hline
Transmission line weights $\Gamma_i$ & \SI{6e-1}{\ohm} \\
Capacitance weight $C_i$, $i=1, 3$ & 
$4 \times 10^{-2} \sqrt{\textrm{kg}} \textrm{m}/\textrm{s}$
\\
\phantom{Capacitance weight $C_i$,} $i=2$ & 
$8 \times 10^{-2} \sqrt{\textrm{kg}} \textrm{m}/\textrm{s}$
\\
Nominal voltage $V^*$ & \SI{48}{\volt} \\
Integral controller weights $D_i$ & \num{1e-4} \\
Load values $-P_l^*$ & \SI{35}{\watt} \\
\hline
\hline
\end{tabular}
\caption{Simulation parameter values.}
\label{t:simparams}
\end{table}

\section{Conclusions}\label{sec.conclusions}

We have proposed controllers for DC microgrids that average power measurement at the sources. The results apply to network preserved model (systems of DAE) of the microgrid in the presence of ZIP loads. Capacitors at the terminals of the grid that model either $\Pi$-models of {lines} or power converter components can be included by means of passivity-based analysis.

Many interesting  new research directions can be taken. The first one is to consider more complex scenarios such as the inclusion of  {dynamical (inductive) lines and loads}. Another one is the  extensions of {the controllers} to network preserved AC microgrids.     Moreover, although the preservation of the geometric mean of the voltages allows for an  estimate of the voltage excursion, no active voltage regulation is present in the proposed scheme. An addition of {voltage controllers} to the power consensus algorithm is an interesting open problem. The power consensus algorithms lead to a new set of power flow equations, whose solvability {still needs to be investigated, e.g., starting from recent advances concerning power flow feasibility and approximations; see \cite{bolognani.tps16,barabanov,JWSP-FD-FB:14c} and references therein. The distributed averaging integral controller \eqref{dapi} discussed in Remark \ref{rem.dapi} enjoys the nice feature of not requiring power measurements and could be an enthralling algorithm to investigate further.
%
Finally, the power consensus algorithms preserves the weighted geometric mean of the voltages and is thus a compelling application for nonlinear consensus schemes \cite{bauso.scl06,cortes.aut2008}. We believe this connection deserves a deeper investigation.}

%

\bibliographystyle{unsrt}
\bibliography{mybiblio_arxiv_v2}

\end{document}